\documentclass[letterpaper, 10 pt, conference]{ieeeconf}

\IEEEoverridecommandlockouts  
\title{\LARGE \bf
Optimal Contract Design for End-of-Life Care Payments
}
\usepackage{amsthm}
\usepackage{graphicx} % Required for inserting images
\usepackage{tikz}
\usepackage{amsmath, amssymb, array, booktabs}
\usepackage{autobreak}
\usepackage{algorithmic}
\usepackage{algorithm}

\usepackage{bbold}
\usepackage{float}

\newtheorem{Theorem}{Theorem}
\newtheorem{Assumption}{Assumption}

\newtheorem{Prop}{Proposition}

\newtheorem{Remark}{Remark}

\author{Muyan Jiang, Ying Chen, Xin Chen, Javad Lavaei, and Anil Aswani% <-this % stops a space
\thanks{*This material is based upon work supported by the National
Science Foundation under Grant DGE-2125913 and Grant
CMMI-1847666.}
\thanks{MJ, YC, XC, JL, and AA are with the Department of Industrial Engineering \& Operations Research, University of California,
Berkeley, CA 94720, USA % <-this % stops a space
{\tt\small \{muyan\_jiang, ying-chen, xin.chen2024, lavaei,aaswani\}@berkeley.edu}}
}

\begin{document}

\maketitle
\thispagestyle{empty}
\pagestyle{empty}

%%%%%%%%%%%%%%%%%%%%%%%%%%%%%%%%%%%%%%%%%%%%%%%%%%%%%%%%%%%%%%%%%%%%%%%%%%%%%%%%
\begin{abstract}
A large fraction of total healthcare expenditure occurs due to end-of-life (EOL) care, which means it is important to study the problem of more carefully incentivizing necessary versus unnecessary EOL care because this has the potential to reduce overall healthcare spending. This paper introduces a principal-agent model that integrates a mixed payment system of fee-for-service and pay-for-performance in order to analyze whether it is possible to better align healthcare provider incentives with patient outcomes and cost-efficiency in EOL care. The primary contributions are to derive optimal contracts for EOL care payments using a principal-agent framework under three separate models for the healthcare provider, where each model considers a different level of risk tolerance for the provider. We derive these optimal contracts by converting the underlying principal-agent models from a bilevel optimization problem into a single-level optimization problem that can be analytically solved. Our results are demonstrated using a simulation where an optimal contract is used to price intracranial pressure monitoring for traumatic brain injuries.
\end{abstract}

%%%%%%%%%%%%%%%%%%%%%%%%%%%%%%%%%%%%%%%%%%%%%%%%%%%%%%%%%%%%%%%%%%%%%%%%%%%%%%%%
\section{INTRODUCTION}

End-of-life (EOL) care is a large part of the total spending on healthcare. For instance, in the United States a significant portion of health expenditures occurs in the final six months of life, comprising approximately 10\% of total national healthcare spending \cite{emanuel1996cost,duncan2019medicare}, with 40\% of this spending occurring in the last 30 days. Despite the enormous spending on EOL care, studies suggest these expenditures do not necessarily improve patient outcomes \cite{iyer2020end,luta2021evidence,anderson2019communication}. Given the significant disconnect between the current level of spending on EOL care and patient outcomes, it is important to study further how to balance these issues better.

Typically, EOL decisions are informed by discussions with key stakeholders, including partners, family members, or the patients themselves. Evidence suggests that early EOL conversations correlate with less aggressive and more cost-effective EOL care, with as much as a 95\%  reduction in expenditure \cite{starr2019associations}, and receipts of palliative care can reduce hospitalization costs and the likelihood of readmission \cite{doi:10.1089/jpm.2013.0612}.
%However, initiating these discussions and implementing the corresponding treatment plans requires action on the part of the healthcare providers. This paper aims to study whether payment schemes can be designed to encourage healthcare providers to carefully consider the medical usefulness and value of providing EOL care.

At a broad level, a number of different payment models have been proposed for healthcare, including capitation, salary, fee-for-service (FFS), pay-for-performance (P4P), and combinations thereof \cite{gosden1996capitation}; because it has been recognized that changes in how healthcare providers are paid can significantly influence the cost-effectiveness of clinical decisions \cite{donaldson1989paying}. In the realm of EOL care, there has been a growing advocacy for integrating pay-for-performance incentives to elevate the quality of palliative care services \cite{bernacki2012improving}. Concurrently, empirical evidence suggests a preference for the fee-for-service model for terminally ill patients, particularly when compared with other Medicare plans \cite{slutsman2002managing}.

\subsection{Related Work on Payment Models}
Current incentives result in more use of services \cite{care2015dying}, and several works have studied contract design for related topics. Research has explored the influence of medical choice on healthcare costs in chronic disease patients, examining the interaction between these choices and principal-agent dynamics \cite{li2022effect}. Additionally, the issue of common-agency problems in the US healthcare system has been scrutinized, particularly in terms of its impact on healthcare contracting and care coordination \cite{frandsen2019sticking}. Furthermore, the design and implications of financial incentives in integrated care systems, with a particular focus on bundled care models, have been analyzed \cite{kadu2021designing}. Overall, advanced medicare payment leads to greater patient quality of life and satisfaction, cost savings, and provider satisfaction \cite{SONENBERG2018112}. Despite these studies, optimal payment pricing for EOL care has not been studied to the best of our knowledge.

\subsection{Payment Models for EOL Care}

The main challenge with designing efficient payment models for EOL care is a fundamental information asymmetry that occurs between the party paying (e.g., health insurance or government) and the healthcare provider. The provider has more knowledge about patients than the payer, so the payer cannot directly observe whether some provided EOL care was necessary or unnecessary. The payer only observes the amount of effort (i.e., the procedures and their costs) that the provider exerts into the EOL care and the outcomes of the patient (e.g., lived or died within 30 days). As a result, it is common for the payer to use fee-for-service as a payment model whereby the provider is paid commensurate with the costs of the procedures performed. 

The setup of the above-described situation motivates the core idea of our paper: Suppose we design a contract where the payment amount to a provider depends upon both the amount of effort and the outcome of the patient. Would such a contract better incentivize providers to only provide necessary EOL care? The idea is that the contract could be designed such that it gives the highest payments when the provider exerts effort and the patient has a good outcome, whereas lower payments are provided when the patient has a poor outcome. By modulating payments for effort by the outcome, our proposed contract combines elements of fee-for-service and pay-for-performance payment models.

\subsection{Contributions and Outline}
Our paper makes three contributions. First, to the best of our knowledge, the idea of modulating EOL care payments based on a binary outcome measure (e.g., survival) is novel. The underlying principal-agent models we propose are also new. Second, we derive the optimal contracts under these models. Though the principal-agent models are described by a bilevel optimization problem, we implicitly convert these into single-level optimization problems by using constraints that capture the fact that healthcare providers maximize their utility functions. This allows us to solve the optimization problem analytically to derive the optimal contract. Third, we propose an algorithm to estimate the parameters appearing in our optimal contracts from real-world data. This estimation methodology is demonstrated using real data.

Sect. II starts by defining the utility functions for the payer and provider, and then it formulates three different principal-agent models, where each model considers a different level of risk tolerance for the provider. In Sect. III, we solve the bilevel optimization problem underlying our principal-agent models in order to derive the optimal contract. We also discuss qualitative insights generated from the optimal contracts. Finally, Section V offers a practical demonstration of our framework through a numerical simulation. This simulation is contextualized in the setting of intracranial pressure monitoring for traumatic brain injuries, and it features a novel algorithm for the estimation of model parameters.

%%%%%%%%%%%%%%%%%%%%%%%%%%%%%%%%%%%%%%%%%%%%%%%%%%%%%%%%%%%%%%%%%%%%%%%%%%%%%%%%
\section{MODEL FORMULATIONS}
In the healthcare system, providers deliver EOL care to patients and receive compensation from payers. Patients are categorized into two groups: $S = 1$ denotes a favorable/good responder status, while $S = 0$ indicates an unfavorable/bad responder status. Identified good responder status indicates patients' profiles are associated with more extensive treatment benefits. There are various statistical models for classifying responder status using the baseline demographic and clinical characteristics and predicting heterogeneous response to a certain treatment applied to a disease \cite{bovis2019defining}. With the knowledge about the patients' conditions, healthcare providers decide expenditure level $E$ regarding whether to administer intensive interventions incurring high costs $E=1$ or opt for palliative care with lower expenditures $E=0$. Consequently, the patient outcomes are reflected in $Q=1$ for survival and $Q=0$ for mortality. The ultimate reimbursement $P = p_{ij}$ is contingent upon both the patient's outcome $Q=i$ and the healthcare expenditure level $E=j$ for $i,j\in \{0,1\}$. Therefore, the healthcare provider must carefully weigh both the quality and cost aspects of end-of-life care.

\subsection{Payer's Utility Function}
We assume a Bernoulli model for the patient's outcome $Q$ depending on responder status $S=s$ and expenditure level $E=j$, $Q\sim \text{Bern}(\pi_{sj})$ for $s,j\in \{0,1\}$. The probability of a good responder status is also assumed as another independent Bernoulli model $S \sim \text{Bern}(\gamma)$ for $\gamma \in (0,1)$.

\begin{Assumption}
$\pi_{01} \geq \pi_{00}, \pi_{11} \geq \pi_{10}, \pi_{10} \geq \pi_{00}, \pi_{11} \geq \pi_{01}$
\label{asmp-pi}
\end{Assumption}
\begin{Remark}
This indicates that higher expenditure generally results in a better survival rate. Moreover, a good responder yields a higher survival rate than a bad responder at the same expenditure level. Additionally, under the Bernoulli assumption, $0<\pi_{01}$, $\pi_{00}$, $\pi_{11}$, $\pi_{10}$, $\gamma<1$.
\end{Remark}

The payer's utility is a weighted sum of the expected survival rate of the patients and the expected payment:
\begin{align}
    u_{\text{{payer}}} = \mathbb{E}(Q) -  \phi\cdot\mathbb{E}(P),
    \label{eq:payer-utility}
\end{align}
where $\mathbb{E}(Q) = (1 - \gamma)\cdot(1 - \pi_{00})+ \gamma\cdot\pi_{11}$ and $\mathbb{E}(\cdot)$ denotes expectation. Here, that constant parameter $\phi>0$ is a weight.
\subsection{Provider's Utility Function}
The provider incurs a disutility $F$ for high expenditure. Without loss of generality, we normalize the units of the principal-agent models by assuming that $F = 1$. Considering the payments and costs, the provider's utility is
\begin{align}
u_{\text{{provider}}} &= \mathbb{E}(P|E) -\mathbb{1}(E=1)\cdot F
% &= \mathbb{1}(E=1)[\gamma\cdotp_{11}+(1-\gamma)\cdotp_{01}-F]\\
% &\quad+\mathbb{1}(E=0)[\gamma\cdotp_{01}+(1-\gamma)\cdotp_{00}]
\label{eq:provider-utility}
\end{align}
where $\mathbb{1}(\cdot)$ is an indicator function that equals one when the condition inside is true and zero otherwise.

\subsection{Principal-Agent Model Formulations}
%%%%%%%%%%%%%%%%%%%%%%%%%%%%%%%%%%%%%%%%%%%%%%%%%%%%%%%%%%%%%%%%%%%%%%%%%%%%%%%%
We construct principal-agent models where the payer is the principal, and the provider is the agent. We constrain our model to ensure the expenditure level aligns with the patient type; that is, the healthcare provider will exert intensive interventions on good responders and palliative care otherwise. 

\subsubsection{Free Payment Model}
The first model is to maximize the payer's utility while allowing the payer to \emph{fine} the provider, if desired at optimality, for bad outcomes:
\begin{align}
\begin{split}
\max_{p_{ij}} \quad & \mathbb{E}(Q)-  \phi\cdot\mathbb{E}(P)\\
\text{s.t.} \quad 
& \mathbb{E}(P)\geq0\\
&  \mathbb{E}(P | S=1, E=1) -F \geq \mathbb{E}(P | S=1, E=0)
% \label{eq:free-constraint-e1}
\\
& \mathbb{E}(P | S=0, E=0) \geq \mathbb{E}(P | S=0, E=1)-F 
% \label{eq:free-constraint-e0}
\end{split}
\label{model:free-payment}
\end{align}
The first constraint in \eqref{model:free-payment} ensures that the healthcare provider does not incur any loss (in expectation) in this scheme.

\subsubsection{Non-Negative Payment Model}
The second model does not allow the payer to fine the provider:
\begin{align}
\begin{split}
\max_{p_{ij}} \quad & \mathbb{E}(Q)-  \phi\cdot\mathbb{E}(P)  \\
\text{s.t.} \quad &  \mathbb{E}(P | S=1, E=1) -F \geq \mathbb{E}(P | S=1, E=0)
\\
& \mathbb{E}(P | S=0, E=0) \geq \mathbb{E}(P | S=0, E=1)-F \\
& p_{00}, p_{01}, p_{10}, p_{11}\geq0
\end{split}
\label{model:non-negative-payment}
\end{align}
We also examine the optimal solution under scenarios where a false diagnosis of responder status could occur.

\subsubsection{Risk-Averse Agent Model}

If the provider is risk-averse, then their utility for a payment is concave:
\begin{align}\label{model:riskaverse}
\begin{split}
\max_{p_{ij}} \quad & \mathbb{E}(Q)-  \phi\cdot\mathbb{E}(P)  \\
\text{s.t.} \quad &  \mathbb{E}(g(P) | S=1, E=1) -F \geq \mathbb{E}(g(P) | S=1, E=0)
\\
& \mathbb{E}(g(P) | S=0, E=0) \geq \mathbb{E}(g(P) | S=0, E=1)-F \\
& p_{00}, p_{01}, p_{10}, p_{11}\geq0
\end{split}
\end{align}
where $g(\cdot)$ is a bijective concave function that is positively valued whenever its argument is a positive value.

\section{OPTIMAL CONTRACTS}
Next, we design optimal contracts for the models described above by solving the corresponding optimization problems.
\subsection{Free Payment Model}
For every feasible incentive design, $\mathbb{E}(Q)$ is a constant with respect to the outcome probabilities. In this situation, we want the optimal solution with the smallest expected reimbursement. The preferred optimal solution can be reduced to solving the following system of equations:
\[
\left\{
\begin{array}{l}
\mathbb{E}(P) = 0 \\
\mathbb{E}(P | S=1, E=1) - 1 \geq \mathbb{E}(P | S=1, E=0) \\
\mathbb{E}(P | S=0, E=0) \geq \mathbb{E}(P | S=0, E=1) - 1
\end{array}
\right.\]
Denote
\begin{align*}
    &\vec{c_0} = [(1 - \gamma)(1 - \pi_{00}), \gamma(1 - \pi_{11}), (1 - \gamma)\pi_{00}, \gamma\pi_{11}]^T\\
    & \vec{c}_1 = [\pi_{10} - 1, 1-\pi_{11}, -\pi_{10}, \pi_{11}]^T\\
    & \vec{c}_2 = [1 - \pi_{00}, \pi_{01}-1, \pi_{00},-\pi_{01}]^T\\
    &b_0 = 0, b_1 = 1,  b_2 = - 1
\end{align*}
With the definition for $\vec{P} = \left[\begin{array}{cccc}p_{00}&p_{01}& p_{10}& p_{11} \end{array}\right]^T$, the problem is equivalent to solving:
\begin{align}\label{model:norm-basic}
\begin{split}
\left\{
\begin{array}{l}
    \vec{c_0} \cdot \vec{P} = b_0\\
    \vec{c_1} \cdot \vec{P} \geq b_1 \\
    \vec{c_2} \cdot \vec{P} \geq b_2
\end{array}
\right.
\end{split}
\end{align}
\begin{Prop}\label{prop:binding}
    The binding system of \eqref{model:norm-basic} admits a solution if and only if $ (1 - \gamma) \pi_{01} + \gamma \pi_{11} < 1$.
    %and $(1 - \gamma) \pi_{00} + \gamma \pi_{10} < 1$.
\end{Prop}

\begin{proof}
    For the system
        \begin{align*}
        \begin{bmatrix}
            \vec{c_0}^T \\
            \vec{c_1}^T \\
            \vec{c_2}^T
        \end{bmatrix}
        \vec{P} =
        \begin{bmatrix}
            b_0 \\
            b_1 \\
            b_2
        \end{bmatrix}
    \end{align*}
    One can do row reduced echelon on the coefficient matrix      $\begin{bmatrix}
        \vec{c_0}^T \\
        \vec{c_1}^T \\
        \vec{c_2}^T
    \end{bmatrix}$ and arrive at:
        \begin{align*}
            \left(
        \begin{array}{cccc}
         1 & 0 & 0 & -\frac{(\pi_{01}-\pi_{11}) ((\gamma-1) \pi_{00}-\gamma \pi_{10})}{(\pi_{00}-\pi_{10}) ((\gamma-1) \pi_{01}-\gamma \pi_{11}+1)} \\
         0 & 1 & 0 & \frac{1}{(\gamma-1) \pi_{01}-\gamma \pi_{11}+1}-1 \\
         0 & 0 & 1 & -\frac{(\pi_{01}-\pi_{11}) ((\gamma-1) \pi_{00}-\gamma \pi_{10}+1)}{(\pi_{00}-\pi_{10}) ((\gamma-1) \pi_{01}-\gamma \pi_{11}+1)} \\
        \end{array}
        \right)
        \end{align*}
    which shows that the matrix has full rank under the condition $ (1 - \gamma) \pi_{01} + \gamma \pi_{11} < 1$.
\end{proof}

\begin{Theorem}
    One optimal solution of \eqref{model:norm-basic} is 
\begin{align*}
    \begin{split}
        p_{00}\to &[-p_{11} (\pi_{11}-\pi_{01}) s_0+\gamma (\pi_{00} \pi_{01}-2 \pi_{00} \pi_{11}\\
        &+\pi_{00}+\pi_{10} \pi_{11}-\pi_{10})+\pi_{00} (\pi_{11}-\pi_{01})]\\
        &/[(\pi_{10}-\pi_{00}) (1-s_1)],
    \end{split}
    \\
    \begin{split}
        p_{01}\to &({-p_{11} s_1-\gamma+1})/({1-s_1}),
    \end{split}
    \\
    \begin{split}
        p_{10}\to &[p_{11} (\pi_{11}-\pi_{01})(1-s_0)+\gamma (\pi_{00} \pi_{01}-2 \pi_{00} \pi_{11}\\
        &+\pi_{00}-\pi_{01}+\pi_{10} \pi_{11}-\pi_{10}+\pi_{11})-(1-\pi_{00}) \\
         &\cdot(\pi_{11}-\pi_{01})]/[{(\pi_{10}-\pi_{00}) (1-s_1)}],
    \end{split}
\end{align*}
where $s_0 = (1 - \gamma) \pi_{00} + \gamma \pi_{10} < 1$ and $s_1 = (1 - \gamma) \pi_{01} + \gamma \pi_{11} < 1$ are expected survival rate for high and low expenditure respectively. Notice they increase as $\gamma$ increases.
\end{Theorem}

\begin{proof}
    After obtaining the reduced row echelon form in the previous proposition, one can derive the linear solution subspace with $p_{11}$ as a free variable.
\end{proof}

\begin{Remark}\label{rm:p11change}
    Regarding $p_{00},p_{10},p_{10}$ as a function of $p_{11}$ we have
    \begin{align*}
        \frac{\partial p_{00}}{\partial p_{11}}& = -\frac{(\pi_{11}-\pi_{01})s_0}{(\pi_{10}-\pi_{00})(1-s_1)} < 0, \\
        \frac{\partial p_{01}}{\partial p_{11}}& = -\frac{s_1}{1-s_1} < 0,\\
        \frac{\partial p_{10}}{\partial p_{11}}& = \frac{(
        \pi_{11}-\pi_{01})(1-s_0)}{(\pi_{10}-\pi_{00})(1-s_1)} > 0
    \end{align*}
    which means that as we increase the welfare of high expenditure spending with desirable outcomes, to maintain optimality, we need to increase that of low expenditure of survival outcome and decrease others.
\end{Remark}

Proposition \ref{prop:binding} shows that an optimal contract ensures zero expected payment for the payer while upholding the constraints that ensure low expenditures for bad responders and high expenditures for good responders. The zero expected payment is achievable in this model because the payer is able to fine the provider when outcomes are poor.

\subsection{Non-Negative Payment Model}
The optimization problem is equivalent to
\begin{equation}\label{model:nonneg}
\begin{aligned}
\min_{\vec{P}} \quad & c_0^T \vec{P} -b_0 
\\
\text{s.t.} \quad & 
\left[\begin{array}{c}c_1^T\\c_2^T\end{array}\right]\vec{P}
\geq\left[\begin{array}{c}b_1 \\b_2 \end{array}\right]
\\
& \vec{P}\geq0
\end{aligned}
\end{equation}
By adding two slack variables $\vec{V} = \left[\begin{array}{cc}v_1 &v_2 \end{array}\right]^T$ to the above linear programming problem, the model is 
\begin{align}
\begin{split}
\min_{\vec{P},\vec{S}} \quad & \left[\begin{array}{ccc}c_0^T&0&0\end{array} \right]
\left[\begin{array}{cc} \vec{P}^T&\vec{V}^T\end{array}\right]^T-b_0
\\
\text{s.t.} \quad & 
\left[\begin{array}{ccc}c_1^T &-1 &0\\c_2^T &0 &-1\end{array}\right]
\left[\begin{array}{cc} \vec{P}^T&\vec{V}^T\end{array}\right]^T
=\left[\begin{array}{c}b_1 \\b_2 \end{array}\right] \\
& \left[\begin{array}{cc} \vec{P}^T&\vec{V}^T\end{array}\right]\geq0
\end{split}
\label{model:non-negative-slack}
\end{align}
\begin{Assumption}
    ${\pi_{01}}{\pi_{10}} \neq {\pi_{00}}{\pi_{11}}$
\label{asmp-alphabeta}
\end{Assumption}
\begin{Remark}
The assumption $\frac{\pi_{10}}{\pi_{00}} \neq \frac{\pi_{11}}{\pi_{01}}$ means that
the transitional benefit in survival rate from good responder to bad responder is different with different expenditure levels.
\end{Remark}
\begin{Theorem}\label{thm:model2opt}
    Under Assumptions \ref{asmp-pi} and \ref{asmp-alphabeta}, the optimal value of \eqref{model:non-negative-slack} is $\gamma$, with solution of the following form:
    \begin{align*}
        &p_{00}\to 0,
        \\&p_{01}\to t,
        \\&p_{10}\to 0,
        \\&p_{11}\to \frac{1}{\pi_{11}} - \frac{1-\pi_{11}}{\pi_{11}}t,
    \end{align*}
    where $0\leq t\leq 1$.
\end{Theorem}
\begin{proof}
    By adding slack variables, the optimal solution of model \eqref{model:non-negative-slack} must satisfy four active linearly independent constraints. After checking all basic points, 
    and under Assumptions 1 and 2, there are only two solutions both primal feasible and dual feasible, which are
        \begin{align*}
        &p_{00}= 0, p_{01}= 1, p_{10}= 0,p_{11}= 1,v_1 = 0, v_2 = 0,\\
        &p_{00}= 0, p_{01}= 0, p_{10}= 0,p_{11} = \frac{1}{\pi_{11}},v_1= 0, v_2 = 1-\frac{\pi_{01}}{\pi_{11}}
    \end{align*}
    The optimal solution lies on the line segment of these two basic points, with optimal value $m^* = \gamma$.
\end{proof}

The optimal contract \eqref{thm:model2opt} for the non-negative payment model is characterized by a parameter $t$, and the quality of the contract is equivalent for any $t\in[0,1]$. Now consider the case $t = 0$. Here, the provider receives zero payment if there is low expenditure or if there is high expenditure but a bad outcome. The provider only receives a payment when there is a high expenditure and a good outcome. 

In practical scenarios, the assumption of perfect classification of responders is not realistic. To address this, we introduce two parameters: \( w_0 := \Pr(S = 0 | TS = 1) \), representing the false negative rate for responder classification, and \( w_1 := \Pr(S = 1 | TS = 0) \), denoting the false positive rate, where \( TS \) is the true responder status class. These parameters encapsulate the probabilistic inaccuracies inherent in the classification process. It is important to note that this information remains concealed from the healthcare provider, who bases decisions solely on the observed responder class and selects the corresponding level of expenditure. Consequently, this modification primarily impacts the objective value, as it alters the survival distribution, thereby affecting the utility of the payer. This adjustment introduces a more realistic and nuanced dimension to the model, acknowledging the uncertainties present in medical classification processes.

Specifically, the new model is the same as \eqref{model:nonneg} except the objective becomes ${c_0^w}^T \vec{P} - b_0$ where 
\begin{align*}
    {c_0^w}^T = [(1-w_0)(1-\gamma)(1-\pi_{00}) &+ w_1 \gamma (1-\pi_{10}),\\
    w_1(1-\gamma)(1-\pi_{01}) &+ (1-w_0) \gamma (1-\pi_{11}),\\
    (1-w_1)(1-\gamma)\pi_{00} &+ w_0 \gamma \pi_{10},\\
    w_1(1-\gamma)\pi_{01} &+ (1-w_0) \gamma \pi_{11}]^T
\end{align*}
The coefficient is formulated based on the unchanged constraint to incentivize healthcare providers to treat observed favorable responders.
\begin{Prop}\label{prop:w01optimal}
    In \eqref{model:nonneg}, replacing $c_0$ with $c_0^w$, under Assumptions 1 and 2, we obtain the optimal solution at
    \begin{align*}
    p_{00}= 0, p_{01}= 0, p_{10}= 0,p_{11} = \frac{1}{\pi_{11}},v_1= 0, v_2 = 1-\frac{\pi_{01}}{\pi_{11}}
    \end{align*}
    with optimal value ${m_w}^* = \gamma \left(1-\frac{\pi_{01} w_1}{\pi_{1}}-w_0\right)+\frac{\pi_{01} w_1}{\pi_{1}}$.
\end{Prop}

\begin{proof}
    Similar to Theorem \ref{thm:model2opt}, checking all basic points will reveal the only optima in this case.
\end{proof}

Interestingly, the optimal contract for this modified model where patients' responder status may be misclassified by the provider is the same as the contract corresponding to $t = 0$ for the original form of this model with exact patient statuses.

\begin{Remark}\label{rm:w01change}
    From Theorem \ref{thm:model2opt} and Proposition \ref{prop:w01optimal}, direct algebra reveals that
    \begin{equation*}
        m^* < m_w^* \Leftrightarrow \frac{\pi_{01}w_1}{\pi_{11}w_0+\pi_{01}w_1} > \gamma
    \end{equation*}
    % This means that the objective value with misclassified responders is higher if and only if 
\end{Remark}

%%%%%%%%%%%%%%%%%%%%%%%%%%%%%%%%%%%%%%%%%%%%%%%%%%%%%%%%%%%%%%%%%%%%%%%%%%%%%%%%

\subsection{Risk-Averse Agent Model}

Considering the model in \eqref{model:riskaverse}, we can define $g^{-1}$ since it is a bijective function. Relabeling $W = g(P)$, and noting the survival rate is again a constant under the constraints, we can reformulate the problem as follows:
\begin{align}\label{model:riskaverse-w}
\begin{split}
\min_{w_{ij}} \quad & \mathbb{E}(g^{-1}(W))  \\
\text{s.t.} \quad &  \mathbb{E}(W | S=1, E=1) -F \geq \mathbb{E}(W | S=1, E=0)
\\
& \mathbb{E}(W| S=0, E=0) \geq \mathbb{E}(W | S=0, E=1)-F \\
& w_{00}, w_{01}, w_{10}, w_{11}\geq0
\end{split}
\end{align}
which is equivalent to
\begin{align}\label{model:riskaverse-vector}
\begin{aligned}
\min_{{W}} \quad & c_0^T g^{-1}({W}) -b_0 
\\
\text{s.t.} \quad & 
\left[\begin{array}{c}c_1^T\\c_2^T\end{array}\right]{W}
\geq\left[\begin{array}{c}b_1 \\b_2 \end{array}\right]
\\
& {W}\geq0
\end{aligned}
\end{align}
Moreover, we can compute solutions satisfying the KKT conditions and claim their optimality because $g^{-1}$ is convex due to the concavity of $g$ and the slater's condition is satisfied by Assumption \ref{asmp-alphabeta} in the constraints.

\begin{Prop}
$p_{00}= g^{-1}(0)$, $p_{01}= g^{-1}(1)$, $p_{10}= g^{-1}(0)$, $p_{11}= g^{-1}(1)$ is an optimal solution of problem \eqref{model:riskaverse-w} with optimal optimal value $\gamma g(1)$.
\end{Prop}

\begin{proof}
    Introduce $\lambda_1,\lambda_2$ as the Lagrangian multipliers of the two inequality constraints, and $\mu = [\mu_{00},\mu_{01},\mu_{10},\mu_{11}]^T$ as the Lagrangian multipliers of the non-negativity constraints with corresponding subscripts.

    We can write out the KKT conditions as follows:
    \begin{equation}
    \left\{
    \begin{aligned}
    &c_0 \cdot \nabla g^{-1}(W) - \lambda_1 c_1 - \lambda_2 c_2 - \mu = 0 \\
    &\lambda_1 (c_1^T W - b_1) = 0 \\
    &\lambda_2 (c_2^T W - b_2) = 0 \\
    &\mu_{ij} W_{ij} = 0 \quad \text{for } i,j = 0,1\\
    & \lambda_1,\lambda_2,\mu \geq 0
    \end{aligned}
    \right.
    \end{equation}
    Any point satisfying this KKT condition is optimal. Specifically, when $W = [0,1,0,1]^T$, we have $c_1^T W - b_1 = c_2^T W - b_2 = 0$, and thus $\mu_{01}=\mu_{11} = 0$, and $\lambda_1,\lambda_2, \mu_{00},\mu_{11} \geq 0$. Denote $c_j[k]$ as the $k^{\text{th}}$ coordinate of vector $c_j$. The Lagrangian becomes:
    \begin{equation}\label{eq:kkt}
    \left[
    \begin{array}{cccc}
    c_{1}[1] & c_{2}[1] & 1 & 0 \\
    c_{1}[2] & c_{2}[2] & 0 & 0 \\
    c_{1}[3] & c_{2}[3] & 0 & 1 \\
    c_{1}[4] & c_{2}[4] & 0 & 0 \\
    \end{array}
    \right]
    \left[
    \begin{array}{c}
    \lambda_1 \\
    \lambda_2 \\
    \mu_{00} \\
    \mu_{01} \\
    \end{array}
    \right]
    =
    \left[
    \begin{array}{c}
    c_0[1] \nabla g^{-1}(0) \\
    c_0[2] \nabla g^{-1}(1) \\
    c_0[3] \nabla g^{-1}(0) \\
    c_0[4] \nabla g^{-1}(1) \\
    \end{array}
    \right]
    \end{equation}
    Solving \eqref{eq:kkt}, we have
\begin{equation}
    \left[
    \begin{array}{c}
    \lambda_1 \\
    \lambda_2 \\
    \mu_{00} \\
    \mu_{01} \\
    \end{array}
    \right]
    =
    \begin{bmatrix}
    \nabla g^{-1}(1) \gamma \\
    0 \\
    \begin{aligned}
    \nabla g^{-1}(0) (1 - \gamma)& (1 - \pi_{00}) \\
    &+ \nabla g^{-1}(1) \gamma (1 - \pi_{10})
    \end{aligned} \\
    \begin{aligned}
    \nabla g^{-1}(0) \gamma& (1 - \pi_{00}) \\
    &+ \nabla g^{-1}(1) \gamma \pi_{10}
    \end{aligned}
    \end{bmatrix}
\end{equation}
which satisfies non-negativity multiplier constraints.
\end{proof}

The main feature of the optimal contract for this model is that payments depend purely upon the expenditure level. Low expenditures get a payment of $g^{-1}(0)$, and high expenditures get a payment of $g^{-1}(1)$. The consequence is that our proposed approach is ineffective in this case because the risk aversion leads to a situation where the providers demand to be compensated when they produce high expenditures, regardless of the outcome. This means it is impossible to incentivize providers to induce low expenditures for bad-responding patients. 

%%%%%%%%%%%%%%%%%%%%%%%%%%%%%%%%%%%%%%%%%%%%%%%%%%%%%%%%%%%%%%%%%%%%%%%%%%%%%%%%

\section{NUMERICAL SIMULATION}

Though intracranial pressure (ICP) monitoring is generally advised for patients with a severe traumatic brain injury (TBI), its impact on patient outcomes is not well-established. There is evidence that ICP monitoring reduces in-hospital and two-week post-injury mortality \cite{aiolfi2017brain}. However, more recent studies question its value because of its cost-effectiveness and risk/benefit ratio \cite{ZAPATAVAZQUEZ201796,mikkonen2020one,anania2021escalation}.

% Critically ill patients may still face mortality despite receiving invasive interventions, including mechanical ventilation. A significant debate exists around whether withholding or withdrawing such treatments might influence the time until death for dying ICU patients \cite{ramazzotti2019withholding}. Additionally, there has been increasing research focus on the economic aspects of mechanical ventilation in the elderly which leads to substantial healthcare costs incurred during the final stages of life \parencite{cox2007economic,jo2023medical,mir2021palliative}.

To numerically evaluate our optimal contracts, we first extracted a cohort of 25934 patients from the MIMIC-IV database (Medical Information Mart for Intensive Care, version 4) \cite{johnson2023mimic}. Out of the cohort, 728 patients were assigned ICP monitoring. The cohort selection criteria are as follows: admitted to ICU due to traumatic brain injury or neurological disease; no relevant data is missing; has available data on the Glasgow Coma Score (GCS). The admission diagnosis is identified by the International Classification of Diseases (ICD) code in the database. For patients with multiple ICU stay records in the database, only the first ICU stay is kept. 

In our framework, we aim to classify the cohort by the treatment given as well as the response level calculated in the data-driven approach illustrated in the pseudo-code below. The responder score construction is inspired by \cite{bovis2019defining}. We first obtained a refined cohort using 1-1 propensity score matching and obtained 1456 patients. For the treatment and control group, we estimated a Cox proportional hazards model and obtained the log hazard ratio coefficients. Assuming a common baseline hazard function, one can take the difference between treatment and control coefficients and obtain a patient-level response score to the treatment. Depending on the practical need, one can choose a cutoff point (here, we select cutoff = 0) if the response scores and cluster the patients. Lastly, we obtain cluster-level outcome rates.

% We obtain our estimates of survival rates conditioning on prognosis and expenditure level by following the algorithm illustrated in the pseudo-code. Since MIMIC-IV is an observational dataset where the decision of mechanical ventilation carries information about prognosis itself, we calculated propensity score to stratified the patients into $K=10$ stratum, based on prognostic scores included in APSIII system (heart rate, mean blood pressure, temperature, respiratory rate, etc.). The propensity scores capture patients' state when admitted and stratification can account for heterogeneity among different groups of subjects. We eventually obtained stratum-level estimates from a Kaplan-Meier survival analysis and output the averages. 

\begin{algorithm}
\caption{Clustering by Treatment and Response and Calculate Outcome Rates}
\begin{algorithmic}[1]
\label{algo:km}
\REQUIRE The sample space $S$ with the following attributes for each patient $i$:
\begin{enumerate}
    \item Baseline covariates, $Z_i \in \mathbb{R}^p$
    \item Propensity score, $p_i \in [0, 1]$
    \item Expenditure level (treatment-control indicator), $e_i \in \{0,1\}$
    \item Death in days, $t_i \in \mathbb{Z}^+$
\end{enumerate}
\STATE 1-1 PS matching based on $p_i$ for every $e_i = 1$ from $S$ and obtain refined cohort $S'$
\STATE For both treatment and control group in $S'$, fit two separate Cox models
\begin{align*}
\log(h_{0}(t)) &= \log(h_{00}(t)) \\&+ (\beta_{10}Z_1 + \beta_{20}Z_2 + \ldots + \beta_{p0}Z_p) \\
\log(h_{1}(t)) &= \log(h_{01}(t)) \\&+ (\beta_{11}Z_1 + \beta_{21}Z_2 + \ldots + \beta_{p1}Z_p)
\end{align*}

\FOR{each $i$ in $S'$}
    \STATE Compute $D(Z_i) = (\beta_{11} - \beta_{10})Z_{i1f} + \ldots + (\beta_{p1} - \beta_{p0})Z_{ip}$
\ENDFOR
\FOR{each combination of $k$ and $(r, e) \in \{0,1\}^2$}
\STATE Extract outcome rates estimate:
    \[
        \hat{\pi}_{r,e} = \frac{\sum_{e_i = e \text{ , } D(Z_i)^+ = r}f(t_i)}{\sum_{i\in S'}\mathbb{1}(e_i = e \text{ and } D(Z_i)^+ = r)}
    \]
    Where $(\cdot)+ = max(\cdot, 0)$ and $f(\cdot)$ is a criteria indicator function that checks outcome (e.g., $f(t_i) = \mathbb{1}(t_i < los_i)$ checks if a patient dies before discharge by length of stay in ICU).
\ENDFOR
\RETURN $\hat{\pi}_{r,e}$ for each group $(r,e)$
\end{algorithmic}
\end{algorithm}

\begin{figure}[thpb]
    \centering
    \includegraphics[width=1.0\linewidth]{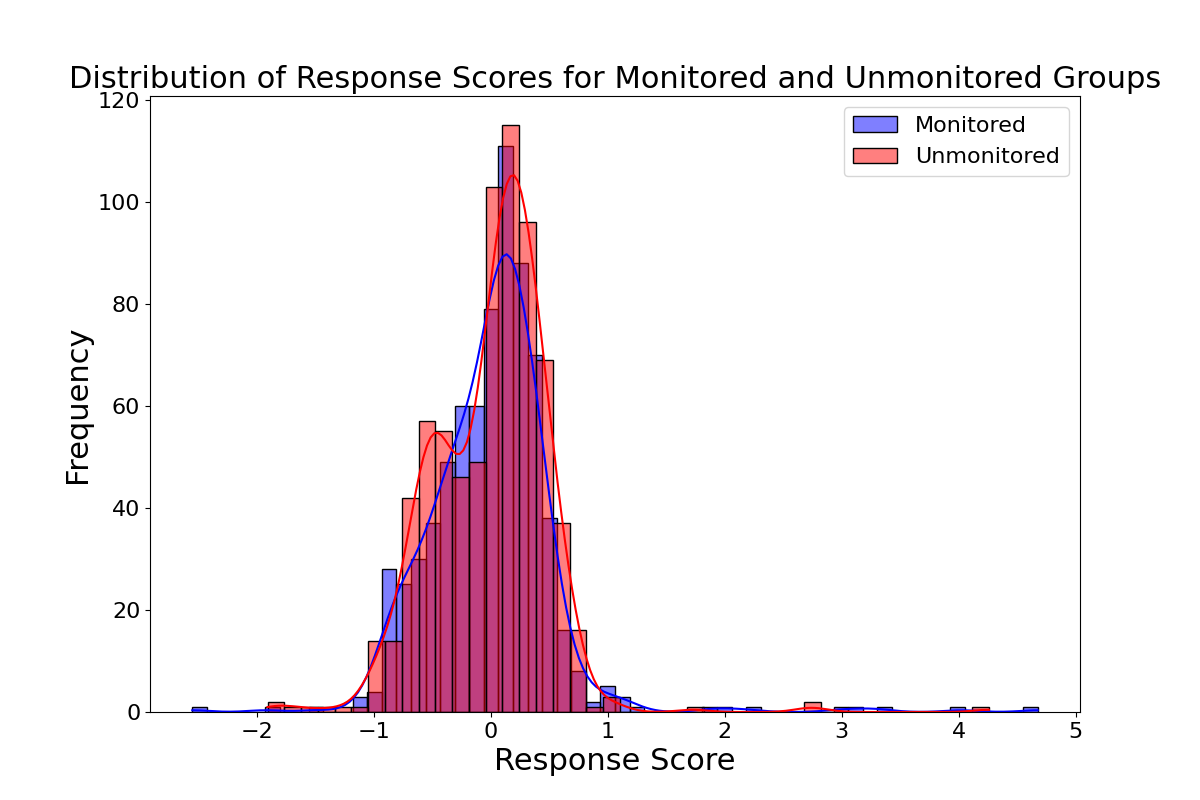}
    \caption{Response Scores Distribution}
    \label{fig:res}
\end{figure}

The response scores estimated from the algorithm are shown in Fig. \ref{fig:res}. The distributions show that the monitoring assignment is not strongly related to the responsiveness of the patients. This lack of reference allows for better expenditure allocation utilizing the information and incentivizes desirable actions. The estimation results of the parameters are shown in Table \ref{tab:outcome_prob}.

\begin{table}[thpb]
\centering
\caption{Outcome Probabilities by Expenditure and Responder Status}
\begin{tabular}{cc}
\hline
$\hat{\pi}$ & Estimates \\
\hline
$\hat{\pi}_{00}$ & $0.51$ \\
$\hat{\pi}_{01}$ & $0.75$ \\
$\hat{\pi}_{10}$ & $0.66$ \\
$\hat{\pi}_{11}$ & $0.85$ \\
\hline
\end{tabular}
\label{tab:outcome_prob}
\end{table}
For estimation of $F$, we cite \cite{berg2004economic}, where ``the average cost fluctuates between €7,600 and €9,000 per hospitalization", and $F= 1$ is the normalization of this cost at around \$10,000. For the estimation of $\gamma$, we estimate it as the proportion of good responder patients with a positive response score in our dataset. $\hat{\gamma} = (335+309)/1456$ = 0.44.

The above estimates satisfy all model parameter assumptions for the non-negative payment model. From the analysis in Sect. IV, if we choose an optimal policy that has the largest incentive gap, the optimal payment is given by $p_{00}= 0, p_{01}= 0, p_{10}= 0,p_{11} = \frac{1}{\hat{\pi}_{11}} = \frac{1}{0.85}\approx 1.18$ with incentive gap $v_2 = 1-\frac{\hat{\pi}_{01}}{\hat{\pi}_{11}} = 1-\frac{0.75}{0.85} = 0.12$ and an objective value $\hat{\gamma} = 0.44$. Recall that with normalized $F = 1$, these numbers should be understood as percentages relative to actual values of $F$, in our case roughly $\$10,000$. Therefore, the payment for high expenditure spending with desirable outcome should be $\$10,000 \times 1.18 = \$11,800$ and obtain a $\$10,000 \times 0.12 = \$1,200$ monetary incentive gap for exerting low expenditure with unfavorable outcomes. This will achieve an optimal expected payment of $\$10,000 \times 0.44 = \$4,400$.

\begin{figure}[thpb]
    \centering
    \includegraphics[width=1.0\linewidth]{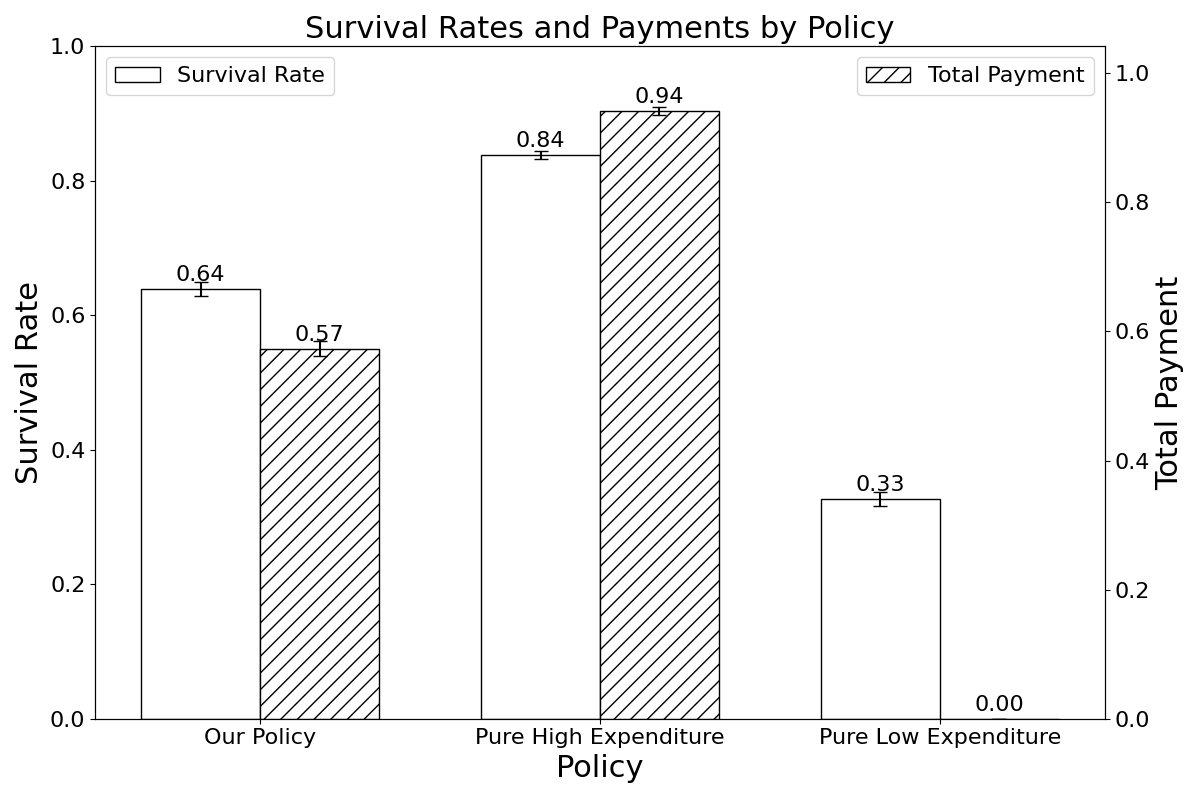}
    \caption{Simluation Comparison with Extreme Policies}
    \label{fig:sim}
\end{figure}

We compare the proposed optimal contract in comparison with two other extreme policies in Fig.\ref{fig:sim} where pure conservative and pure aggressive treatment practices are implemented. 

In the optimal payment plan, where only patients diagnosed with a positive response score get treated, it achieves a 64\% survival rate with a 0.55 payment investment. Compared to ``pure high expenditure" with 83\% survival rate with 0.94 payment, our policy achieves both higher average outcome-rate/payment ratio $(\frac{0.64}{0.55} > \frac{0.83}{0.94})$ as well as marginal outcome-rate/payment ratio when benchmarked against ``pure low expenditure" policy $(\frac{0.64-0.35}{0.55} > \frac{0.83-0.35}{0.94})$, which proves its better cost-effectiveness.

\section{CONCLUSION}
Optimal contracts have been designed across three different models to ensure that payment structures incentivize providers to increase expenditure for good responders while encouraging reduced expenditure for poor responders. Numerical simulation utilizing MIMIC-IV data on ICP monitoring for TBI patients was used to evaluate a new policy. When healthcare providers are generally risk averse, our results show that in this case, the optimal contract is such that the risk aversion leads to a situation where the providers demand to be compensated when they produce high expenditures, regardless of the outcome. This means it is impossible to incentivize providers to induce low expenditures for bad-responding patients. 

\section*{ACKNOWLEDGEMENTS}
The authors thank Malini Mahendra, MD, in the Department of Pediatrics, Division of Pediatric Critical Care, UCSF Benioff Children's Hospital and Philip R. Lee Institute for Health Policy Studies, University of California, San Francisco, CA, USA, for her input about ICP monitoring.
%%%%%%%%%%%%%%%%%%%%%%%%%%%%%%%%%%%%%%%%%%%%%%%%%%%%%%%%%%%%%%%%%%%%%%%%%%%%%%%%

\bibliographystyle{IEEEtran}
\bibliography{ref}

\end{document}